\documentstyle[12pt]{article}
\newcommand{\copyleft}{
GNU FDL\thanks{
Copyright (C) 1990, 1998 Peter G. Doyle.
Permission is granted to copy, distribute and/or modify this document
under the terms of the GNU Free Documentation License, 
as published by the Free Software Foundation;
with no Invariant Sections, no Front-Cover Texts, and no Back-Cover Texts.
}}
\title{The knee-jerk mapping}
\author{
Peter G. Doyle \and Jim Reeds
}
\date{Version dated 5 October 1998
\\ \copyleft}

\newcommand{\xbar}{\bar{x}}
\newcommand{\ubar}{\bar{u}}
\newcommand{\vbar}{\bar{v}}
\newcommand{\xx}{{x_1,\ldots,x_n}}
\newcommand{\xxbar}{{\xbar_1,\ldots,\xbar_n}}
\newcommand{\Zbar}{\bar{Z}}
\newcommand{\bx}{{\bf x}}
\newcommand{\bxbar}{{\bf \xbar}}
\newcommand{\by}{{\bf y}}
\newcommand{\bu}{{\bf u}}
\newcommand{\ba}{{\bf a}}
\newcommand{\qed}{\spadesuit}

\begin{document}

\maketitle

\begin{abstract}
We claim to give the definitive theory of what we call the
`knee-jerk mapping',
which is the basis for a class of optimization algorithms
introduced by Baum, and promoted by Dempster, Laird, and Rubin
under the name `EM algorithm'.
\end{abstract}

\section*{Introduction}

We give the definitive theory of the knee-jerk mapping,
to be defined below.
This mapping has been investigated by many people, most notably Baum
(\cite{baumEagon},
\cite{baumPetrie},
\cite{baumSell},
\cite{baumEtAl},
\cite{baum:inequality}
).

We begin with an example,
taken from
\cite{dempsterEtAl:jerk}.
Suppose you want to locate the maximum of the function
\[
Z(x,y) =  x^{34}y^{38}(1+2x)^{125}
\]
on the {\em 1-simplex} 
(a fancy name for a line segment)
\[
\Sigma
=
\{
x,y>0;\;x+y=1
\}
.
\]
One way you can find it is by iterating the
{\em knee-jerk mapping}
\[
(x,y)
\mapsto
\frac{1}{xZ_x+yZ_y}(xZ_x,yZ_y)
=
\ldots
.
\]
This maps the simplex $\Sigma$ to itself,
and what is notable about the mapping is that it
increases the value of the objective function $Z$.

The {\em one true explanation}
of this ratcheting property of the knee-jerk map,
the explanation that {\em lays bare}
once and for all what is going on here,
is as follows:
Like any polynomial with only positive coefficients,
the function
$Z$ is {\em log-log-convex};
that is, $\log Z$ is convex as a function of $(\log x ,\log y)$;
that is, 
\[
W(u,v) = \log Z(e^u,e^v)
\]
is convex as a function of $(u,v)$.
We're trying to find the maximum of $W$ on the set
\[
T =
\{e^u + e^v = 1\}
.
\]
Since $W$ is convex,
if we fix a point $(u,v)$,
the graph of $W$ lies above its tangent plane at
$(u,v,W(u,v))$:
\[
W(\ubar,\vbar)
\geq
W_u(u,v)(\ubar-u) + W_v(u,v) (\vbar-v)
.
\]
Now ideally we'd like to move from $(u,v)$ 
directly to the point of $T$ where $W(u,v)$ is greatest.
What the knee-jerk mapping does is move instead to the point
where the lower bound on the right hand side of the inequality
above is maximized.
This can't help increasing the objective function, right?

One remarkable fact should be pointed out, though it won't be
gone into below:
While the function $Z$ is log-log-convex, it is nevertheless
{\em log-concave};
that is,
$\log Z$ is concave as a function of $(x,y)$.
(This is true because $Z$ is a product of homogeneous linear functions
with positive coefficients.)
Because $Z$ is log-concave, it has a unique maximum on the
simplex $\Sigma$.
While all polynomials with positive coefficients are log-log-convex,
only very special polynomials are simultaneously log-concave.

A class of log-concave examples fundamentally more exciting than products of linear
functions can be obtained as follows:
Take a connected graph $G$, think of its edges as variables, form for each
spanning tree of $G$ a monomial (of degree one smaller than the number
of vertices of $G$), and form a polynomial $D_G$---the
{\em discriminant} of $G$---by adding up the monomials corresponding
to all spanning trees of $G$.
For example,
if $G$ is a triangle with edges $x,y,z$,
\[
D_G(x,y,z) = x y + x z + y z
.
\]
Discriminants of graphs are always log-concave.
(If you know what a matroid is,
let me add that the discriminant of a regular matroid is log-concave,
but I don't know if the discriminant of a general matroid always is;
my guess is that it isn't.)

Discriminants of graphs are particular cases of the
{\em diagonal discriminants} of Bott and Duffin;
these are always log-concave
(because the determinant function is log-concave when restricted to the
set of positive-definite matrices)
as well as being log-log-convex
(because they are polynomials with positive coefficients).

\section*{Knee-jerk functions}

In real $n$-space, we will denote the positive orthant by $\Pi$
and the closed standard simplex by $\Sigma$:
\[
\Pi = \{\xx>0\}
,
\]
\[
\Sigma = \{\xx>0;\;x_1+\ldots+x_n=1\}
.
\]
We denote their closures by $\bar{\Pi}$ (the non-negative orthant)
and $\bar{\Sigma}$ (the closed standard simplex).

We say that a function $Z(\xx)$ 
from $\Pi$ to the positive real numbers is
{\em log-log-convex}
if $\log Z$ is a convex function of
$u_1 = \log x_1,\ldots,u_n = \log x_n$.
The name comes from the fact that in the case $n=1$
a log-log-convex function is one whose graph appears
convex when drawn on log-log graph paper.
We say that $Z$ is a {\em knee-jerk function}
if $Z$ is increasing (which we take to mean
what some would call `non-decreasing')
and log-log-convex.
For pedantry's sake we require in addition that $Z$ be smooth, and extend
continuously to $\bar{\Pi}$.

\section*{Properties and examples.}

There are many characterizations of convex functions,
but for our purposes the most important is that
a function is convex if and only if its graph lies above all of its
tangent planes.
Thus a smooth function $Z$ is log-log-convex if and only if
for any two points $\bx=(\xx)$ and $\bxbar = (\xxbar)$,
\begin{eqnarray*}
\log \Zbar - \log Z
&\geq&
(\log Z)_{u_1} (\ubar_1 - u_1) + \ldots + (\log Z)_{u_n} (\ubar_n - u_n)
\\&=&
\frac{x_1 Z_{x_1}}{Z} \log \frac{\xbar_1}{x_1}
+ \ldots +
\frac{x_n Z_{x_n}}{Z} \log \frac{\xbar_n}{x_n}
,
\end{eqnarray*}
where $\Zbar = Z(\bxbar)$ and $(\log Z)_{u_1}$ denotes the derivative of 
$\log Z$ with respect to $u_1$, etc.

Using this characterization of log-log-convexity and Jensen's inequality---which
states that for a concave function like $\log$ the weighted average of the
values is littler than the value of the weighted average---we get a proof
that the function
$Z(\xx) = x_1+ \ldots + x_n$ is log-log-convex, and hence
a knee-jerk function:
\begin{eqnarray*}
&&
\frac{x_1 Z_{x_1}}{Z} \log \frac{\xbar_1}{x_1}
+ \ldots +
\frac{x_n Z_{x_n}}{Z} \log \frac{\xbar_n}{x_n}
\\&=&
\frac{x_1}{x_1+\ldots+x_n} \log \frac{\xbar_1}{x_1}
+ \ldots +
\frac{x_n}{x_1+\ldots+x_n} \log \frac{\xbar_n}{x_n}
\\&\leq&
\log \left (
\frac{x_1}{x_1+\ldots+x_n} \frac{\xbar_1}{x_1}
+ \ldots +
\frac{x_n}{x_1+\ldots+x_n} \frac{\xbar_n}{x_n}
\right )
\\&=&
\log \frac{\xbar_1 + \ldots + \xbar_n}{x_1+\ldots+x_n}
\\&=&
\log \frac{\Zbar}{Z}
.
\end{eqnarray*}

Once we know that $x_1+\ldots+x_n$ is a knee-jerk function,
we can easily produce a wealth of other examples by observing that
the class of knee-jerk functions is closed under a variety of operations.
The coordinate functions $\xx$ are knee-jerk functions,
as is any positive constant function.
Products, positive scalar multiples, and positive (possibly fractional) powers
of knee-jerk functions are knee-jerk functions.
So is the composition $Z(Z_1,\ldots,Z_k)$ of a knee-jerk function
$Z(x_1,\ldots,x_k)$
with knee-jerk functions
$Z_1(\xx),\ldots,Z_k(\xx)$,
because the composition of increasing convex functions is increasing and
convex.
And since $x_1+\ldots+x_n$ is a knee-jerk function,
it follows that sums of knee-jerk functions are knee-jerk functions.
Thus any non-zero polynomial with non-negative coefficients is a knee-jerk
function.

\section*{The knee-jerk mapping}

If $Z(\xx)$ is a knee-jerk function,
we define the {\em knee-jerk mapping}
\begin{eqnarray*}
T_Z(\bx)
&=&
\frac{1}{x_1 Z_{x_1} + \ldots + x_n Z_{x_n}}
(x_1 Z_{x_1},\ldots,x_n Z_{x_n})
\\&=&
\frac{Z}{x_1 Z_{x_1} + \ldots + x_n Z_{x_n}}
(x_1 (\log Z)_{x_1},\ldots,x_n (\log Z)_{x_n})
\\&=&
\frac{1}{(\log Z)_{x_1} + \ldots + (\log Z)_{x_n}}
((\log Z)_{u_1}, \ldots, (\log Z)_{u_n})
.
\end{eqnarray*}
(If $Z_{x_1} = \ldots = Z_{x_n} = 0$, we
define $T_Z(\xx) = \frac{1}{x_1+ \ldots + x_n}(\xx)$---or just pretend we
didn't notice.)
Note that when $Z$ is homogeneous of (possibly fractional) degree $d$,
Euler's identity
\[
x_1 Z_{x_1} + \ldots x_n Z_{x_n} = d Z
\]
implies that
\[
T_Z(\bx)
=
\frac{1}{dZ}
(x_1 Z_{x_1},\ldots,x_n Z_{x_n})
.
\]

$T_Z$ maps the positive orthant $\Pi$ to the closed simplex $\bar{\Sigma}$,
and thus restricts to a mapping of $\Sigma$ to
$\bar{\Sigma}$.
It is easy to see that
a point $\bx \in \Sigma$ is fixed by $T_Z$
if and only if it is a critical point of $Z$ on $\Sigma$.
The great thing about the knee-jerk mapping is that if
$\bx$ is not a critical point of $Z$ on $\Sigma$ then
$Z(T_Z(\bx)) > Z$;
this will be proven in the next section.
This makes the knee-jerk mapping
a natural to iterate if you are interested in finding the
maximum of $Z$ on $\Sigma$.
The name `knee-jerk' is partly meant to suggest the automatic way in which
the mapping increases the objective function $Z$.

\section*{The knee-jerk inequality}

Write
\[
\bx' = T_Z(\bx)
\]
and
\[
Z' = Z(\bx`)
.
\]

\proclaim The knee-jerk inequality.
\[
\log \frac{Z'}{Z}
\geq
\frac{x_1 Z_{x_1} + \ldots + x_n Z_{x_n}}{Z}
\left (
x'_1 \log \frac{x'_1}{x_1}
+ \ldots +
x'_n \log \frac{x'_n}{x_n}
\right )
.
\]

{\bf Proof.}
From the characterization of log-log-convexity above,
we have
\[
\log \Zbar - \log Z
\geq
\frac{x_1 Z_{x_1}}{Z} \log \frac{\xbar_1}{x_1}
+ \ldots +
\frac{x_n Z_{x_n}}{Z} \log \frac{\xbar_n}{x_n}
.
\]
Substituting $\bxbar = \bx`$ yields
the knee-jerk inequality.
$\qed$

Recall
(if you don't already know)
that for probability vectors $\bx \in \Sigma, \by \in \bar{\Sigma}$
the {\em I-divergence}
$I(\by;\bx)$ is defined to be
\[
I(\by;\bx) =
y_1 \log \frac{y_1}{x_1}
+ \ldots +
y_n \log \frac{y_n}{x_n}
.
\]
This quantity is always $\geq 0$, 
with equality if and only if $\bx = \by$.
(This follows from an application of Jensen's inequality
similar to that used above to show that $x_1 + \ldots + x_n$
is a knee-jerk function.)

\proclaim Corollary.
If $\bx \in \Sigma$ then
\[
\log \frac{Z'}{Z}
\geq
\frac{x_1 Z_{x_1} + \ldots + x_n Z_{x_n}}{Z}
I(\bx';\bx)
\geq 0
.
\]
In particular,
$Z' > Z$
unless the point $\bx$ is fixed by $T_Z$,
which happens if and only if $\bx$ is a critical point of $Z$ on $\Sigma$.
$\qed$

\section*{What is going on here?}
Say our goal is to maximize $Z$ over $\Sigma$.
We're sitting at some point $\bx$,
and we want to pick a new point $\bxbar \in \bar{\Sigma}$
so as to increase the objective function $Z$ as much as possible.
Since $Z$ is log-log-convex we know that
\[
\log \Zbar - \log Z
\geq
(\log Z)_{u_1} (\ubar_1 - u_1) + \ldots + (\log Z)_{u_n} (\ubar_n - u_n)
.
\]
The knee-jerk idea is to choose $\bxbar \in \bar{\Sigma}$ so as
to make the lower bound on the right of this inequality as large as possible.
That is, we want to do as well as possible using only the value of $Z$
and its derivatives at $\bx$ and the knowledge that $Z$ is a knee-jerk
function.
So we want to choose $\bxbar$ so as to maximize
\[
F(\ubar_1,\ldots,\ubar_n) \equiv
(\log Z)_{u_1} (\ubar_1 - u_1) + \ldots + (\log Z)_{u_n} (\ubar_n - u_n)
\]
subject to the constraint
\[
G(\ubar_1,\ldots,\ubar_n) \equiv
e^{\ubar_1} + \ldots + e^{\ubar_n} = 1
.
\]
The maximum occurs where
\[
\nabla_{\bar{\bu}} G
=
(\xxbar)
\]
is proportional to
\[
\nabla_{\bar{\bu}} F
=
((\log Z)_{u_1}, \ldots, (\log Z)_{u_n})
,
\]
that is, where
\[
\bxbar = T_Z(\bx)
.
\]

{\bf Ruminations.}
When $\bx \in \Sigma$, the fact that $\bx' = T_Z(\bx)$ maximizes
the lower bound for $\Zbar$ implies right away that $Z' \geq Z$,
independently of the hocus-pocus with the I-divergence.
Indeed, the positivity of the I-divergence can now be seen as
a consequence of the fact that $x_1+\ldots+x_n$ is a knee-jerk function.
This is not so surprising, perhaps, since both facts followed from
very similar applications of Jensen's inequality.
But now it appears that $x_1+\ldots+x_n$ is somehow the most important
of all knee-jerk functions.
And why should it be so distinguished?
Because it crops up in the definition of the simplex $\Sigma$.

\section*{Generalizations}

Given
$ \ba = (a_1,\ldots,a_n)$,
$a_1,\ldots,a_n > 0$,
define
\[
\Sigma_{\ba}
=
\{\xx>0;\; a_1 x_1 + \ldots a_n x_n = 1 \}
\]
and define
\[
T_{Z,\ba} : \Pi \rightarrow \bar{\Sigma}_{\ba}
,
\]
\[
\bx'
=
T_{Z,\ba}(\bx)
=
\frac{1}{x_1 Z_{x_1} + \ldots + x_n Z_{x_n}}
(\frac{x_1 Z_{x_1}}{a_1},\ldots,\frac{x_n Z_{x_n}}{a_n})
.
\]
Then the knee-jerk inequality becomes
\[
\log \frac{Z'}{Z}
\geq
\frac{x_1 Z_{x_1} + \ldots + x_n Z_{x_n}}{Z}
\left (
a_1 x'_1 \log \frac{x'_1}{x_1}
+ \ldots +
a_n x'_n \log \frac{x'_n}{x_n}
\right )
.
\]
When
$\bx \in \Sigma$ this becomes
\[
\log \frac{Z'}{Z}
\geq
\frac{x_1 Z_{x_1} + \ldots + x_n Z_{x_n}}{Z}
I((a_1 x'_1,\ldots,a_n x'_n);(a_1 x_1,\ldots,a_n x_n))
\geq 0
.
\]

More interesting,
we can replace the simplex $\Sigma$ with a product of simplices:
Let
\[
Z = Z(x_{1,1},\ldots,x_{1,n_1},\ldots,x_{k,1},\ldots,x_{k,n_k})
.
\]
Let
\[
T = \{ x_{i,j} > 0 ; \; \sum_j x_{i,j} = 1 \}
,
\]
and define
\[
T_Z : \Pi \mapsto \bar{T}
\]
by
\[
x'_{i,j}
=
\frac{x_{i,j} Z_{x_{i,j}}}
{\sum_j x_{i,j} Z_{x_{i,j}}}
.
\]
Then
\[
\log \frac{Z'}{Z}
\geq
\sum_i
\frac{\sum_j x_{i,j} Z_{x_{i,j}}}{Z}
\left (
\sum_j x'_{i,j} \log \frac{x'_{i,j}}{x_{i,j}}
\right )
,
\]
and when $\bx \in T$,
\[
\log \frac{Z'}{Z}
\geq
\sum_i
\frac{\sum_j x_{i,j} Z_{x_{i,j}}}{Z}
I(
(x'_{i,1},\ldots,x'_{i,n_i});
(x_{i,1},\ldots,x_{i,n_i})
)
\geq
0
.
\]

\bibliography{jerk}

\begin{thebibliography}{1}

\bibitem{baum:inequality}
L.~E. Baum.
\newblock An inequality and associated maximization technique in statistical
  estimation for probabilistic functions of {M}arkov processes.
\newblock In {\em Inequalities, Vol. 3}, pages 1--8. Academic Press, New York,
  1972.

\bibitem{baumEagon}
L.~E. Baum and J.~A. Eagon.
\newblock An inequality with applications to statistical estimation for
  probabilistic functions of {M}arkov processes and to a model for ecology.
\newblock {\em Bull. Amer. Math. Soc.}, 73:360--363, 1967.

\bibitem{baumPetrie}
L.~E. Baum and T.~Petrie.
\newblock Statistical inference for probabilistic functions of finite state
  {M}arkov chains.
\newblock {\em Ann. Math. Stat.}, 37:1554--1563, 1966.

\bibitem{baumEtAl}
L.~E. Baum, T.~Petrie, G.~Soules, and N.~Weiss.
\newblock A maximization technique occurring in the statistical analysis of
  probabilistic functions of {M}arkov chains.
\newblock {\em Ann. Math. Stat.}, 41:164--171, 1970.

\bibitem{baumSell}
L.~E. Baum and G.~R. Sell.
\newblock Growth transformations for functions on manifolds.
\newblock {\em Pacific J. Math.}, 27:211--227, 1968.

\bibitem{dempsterEtAl:jerk}
A.~P. Dempster, N.~M. Laird, and D.~B. Rubin.
\newblock Maximum likelihood from incomplete data via the {EM} algorithm.
\newblock {\em Ann. Math. Stat.}, 41:164--171, 1970.

\end{thebibliography}
\bibliographystyle{plain}

\end{document}